\documentclass[10pt]{article}

\usepackage[top=1.3in, bottom=1.3in, left=1.5in, right=1.5in]{geometry}

\usepackage{amsmath, amssymb, amsfonts, amsthm, amscd, mathtools, graphicx, mathtext,color, wrapfig, verbatim, slashed,url}
\usepackage{yfonts}

\usepackage{mathrsfs}
\usepackage{pifont}

\usepackage{graphicx}
\usepackage[font={small,sl}]{caption}
\usepackage[all]{xy}

\usepackage{tikz}
\usetikzlibrary{matrix,arrows}

\usepackage{lipsum}

\let\OLDthebibliography\thebibliography
\renewcommand\thebibliography[1]{
  \OLDthebibliography{#1}
  \setlength{\parskip}{0pt}
  \setlength{\itemsep}{0pt plus 0.0ex}
}
\makeatletter

\@addtoreset{equation}{section}
\makeatother

\newcommand{\beq}{\begin{equation}}
\newcommand{\eeq}{\end{equation}}

\newcommand{\bP}{\ensuremath{\mathbb{P}}}

\newcommand{\rd}{/\!\!/\!\!/\!\!/}

\newcommand{\scQ}{\ensuremath{\mathcal{Q}}}

\newcommand{\fg}{\mathfrak{g}}
\newcommand{\fh}{\mathfrak{h}}

\newtheorem{theorem}{Theorem}[section]

\DeclareMathOperator{\Aut}{Aut}

\newcommand{\Ct}{\mathbb{C}^\times}

\newcommand{\cO}{\mathscr{O}}

\newcommand{\C}{\mathbb{C}}

\newcommand{\bT}{\mathsf{T}}

\newcommand{\nc}{\newcommand}
\nc{\mc}{\mathcal}
\nc{\wh}{\widehat}
\nc{\ghat}{\wh\fg}
\nc{\on}{\operatorname}

\newcommand{\Lfgh}{\wh{^L\fg}}
\nc{\n}{{\mathfrak n}}
\nc{\pone}{{\mathbb C}{\mathbb P}^1}
\nc{\la}{\lambda}
\nc{\al}{\alpha}

\newtheorem{Theorem}{Theorem}

\theoremstyle{definition}
\newtheorem{remark}[theorem]{Remark}
\begin{document}
\baselineskip=28pt  
\baselineskip 0.7cm

\begin{titlepage}


\renewcommand{\thefootnote}{\fnsymbol{footnote}}

\begin{center}

\vspace*{10mm}
\vskip 1in
{\LARGE \bf
Two-parameter Langlands Correspondence
}
\vspace*{10mm}

{\large
Mina Aganagic$^{1,2}$, Edward Frenkel$^{2}$, Andrei Okounkov$^{3}$
}
\\
\medskip

\vskip 0.5cm

{\it
$^1$Leinweber Institute for Theoretical Physics, University of California, Berkeley\\
$^2$Department of Mathematics, University of California, Berkeley\\
$^3$Department of Mathematics, Columbia University 
}

\end{center}

\centerline{{\bf Abstract}}
\medskip

In \cite{AFO} we established, for every simply-laced Lie algebra
$\fg$, a canonical isomorphism between the spaces of deformed
conformal blocks of the deformed ${\mc W}$-algebra ${\mc
  W}_{q,t}(\fg)$ and the quantum affine algebra of $\wh{\fg}$, which
we view as a $q$-deformation of the quantum Langlands
correspondence. This was done by realizing the deformed conformal
blocks of these algebras via the quantum K-theory of the Nakajima
quiver varieties. We also linked this isomorphism to a duality
emerging from the 6d little string theory. Here, we give a brief
survey of these results and propose an extension to the non-simply
laced case, which exhibits a Langlands-type duality.\footnote{To be
published in the Proceedings of the 2025 International Congress of
Basic Science on the occasion of the Frontiers of Science Award to our
paper \cite{AFO}.}

\noindent\end{titlepage}
\setcounter{page}{1} 

\setcounter{tocdepth}{2}


\section{Introduction}\label{sec:zero}

The Langlands Program was launched by Robert Langlands in late 1960s
with the goal of relating certain questions in number theory and
harmonic analysis. It was subsequently generalized to the geometric
setting of complex algebraic curves, commonly referred to as the
geometric Langlands correspondence. Recently, a sheaf-theoretic, {\em
  categorical Langlands correspondence} has been proved in
\cite{Gaitsgory}, following the seminal paper \cite{BD}. In addition,
a function-theoretic, {\em analytic Langlands
  correspondence}, has been introduced in \cite{EFK}.

In our paper \cite{AFO}, we proposed a two-parameter deformation of a
specific manifestation of the geometric Langlands correspondence;
namely, an isomorphism between two spaces of conformal blocks. We
proved it in the simply-laced case. In this paper, we give a survey of
these results and discuss a generalization to the non-simply laced
case.

\setcounter{footnote}{0}
Let $\fg$ be a simple Lie algebra over $\C$ and $^L\fg$ its Langlands
dual.\footnote{We follow the notation of
\cite{AFO}, which is opposite to what has become the standard notation
in the subject.} In the geometric Langlands
correspondence, a crucial role is played by the center of the affine
Kac--Moody algebra $\wh{^L\fg}$ at the critical level $^Lk=-^Lh^\vee$
and its isomorphism with the classical ${\mc W}$-algebra ${\cal
  W}_\infty({\fg})$, established in \cite{FF} (see also
\cite{F:book}).

The level $^Lk$ of $\wh{^L\fg}$ may be deformed away from the critical
value, and at the same time ${\cal W}_\infty({\fg})$ may be deformed
to the quantum ${\mc W}$-algebra ${\cal W}_\beta({\fg})$.

The algebra ${\cal W}_{\beta}({\fg})$ is obtained by the
quantum Drinfeld--Sokolov reduction \cite{FF:ds,BO,FF} of the affine
algebra $\ghat$ at level $k$, where $\beta = m(k+h^\vee)$, $m$ being
the lacing number of $\fg$.\footnote{${\cal W}_{\beta}(\fg)$ agrees
with the notation of \cite{FR}, and it is ${\cal W}_{k}(\fg)$ of
\cite{FB,Frenkel}, where $\beta = m(k+h^\vee)$.}

Hence one is
naturally led to a {\em quantum deformation} of the geometric
Langlands correspondence. Many interesting structures have been
studied in this framework (see the references in \cite{AFO}). In
addition, it has been linked in \cite{KW} to $S$-duality of maximally
supersymmetric 4d gauge theories with gauge groups being the compact
forms of Langlands dual groups $G$ and $^LG$ with the Lie algebras
$\fg$ and $^L\fg$, respectively.

Let's focus on the following feature of the quantum geometric
Langlands correspondence, an isomorphism between the spaces of
conformal blocks of certain modules of two chiral algebras:
\beq\label{geomL} \widehat{^L\fg}_{{}^Lk}\on{-blocks} \qquad
\longleftrightarrow \qquad {\cal W}_{\beta}({\fg})\on{-blocks}, \eeq if the
parameters are related by the formula (here $^Lh^\vee$ is the dual Coxeter
number of $^L\fg$) \beq\label{zero2}
\beta = \frac{1}{^Lk+{}^Lh^\vee} \eeq and are {\em generic}, meaning that
$\beta$ is {\em not} a rational number. (For an explanation of how this
isomorphism of the
spaces of conformal blocks follows from the quantum geometric
Langlands correspondence, see Section 6.4 of \cite{AFO}.)

In \cite{AFO}, we proved a stronger
result in the case of simply-laced $\fg$ and genus zero curve ${\mc
  C}$: a {\em canonical isomorphism} between the spaces of conformal
blocks of the two algebras if we modify the relation \eqref{zero2} to
\beq\label{zero1} \beta - m =
\frac{1}{^L(k+h^\vee)} \eeq
The shift by $-m$ in formula \eqref{zero2} is essential (see
\cite{AFO}, Section 6.4, for more details). In the simply-laced case
$m=1$ and $^Lh^\vee=h$ but we will need formula \eqref{zero1} in the
non-simply laced case (see Section \ref{nonsimplylaced}).

\begin{remark}
The relation between the corresponding chiral algebras and their
conformal blocks may be viewed as a strong/weak coupling
transformation well-known in quantum physics. Indeed, if we define
$\tau = \beta/m$ and $^L\tau = - {^L(k+h^\vee)}$, then \eqref{zero1}
becomes \beq\label{firstl} \tau -1 = -1/(m ^L\tau), \eeq and so small
values of $^L\tau$ correspond to large values of $\tau$. The
parameters $\tau$ and $^L\tau$ are related to the complexified
coupling constants of the two $S$-dual 4d Yang-Mills theories. The
shift of $\beta$ by $-m$ (passing from \eqref{zero2} to \eqref{zero1})
corresponds to the shift $\tau \mapsto \tau-1$ is a shift
of the theta-angle (see \cite{AFO} for more details).\qed
\end{remark}

The canonical isomorphism between the spaces of conformal blocks
established in \cite{AFO} arises only after we introduce one more
parameter and perform one more deformation, which was the main idea of
\cite{AFO}. Namely, we introduced there a {\em two-parameter
  deformation} of the geometric Langlands correspondence, a
$q$-deformation in addition to the deformation away from the
critical level.

To define this two-parameter deformation, we replace the above
algebras with their deformed counterparts. The first is the quantum
affine algebra of level $^Lk$, which we denote by
$U_{\hbar}(\wh{^L\fg})_{^Lk}$.
The second algebra is the deformed ${\mc W}$-algebra ${\cal
  W}_{q,t}({\fg})$ introduced in \cite{FR}, which is a deformation of
the conformal ${\mc W}$-algebra ${\mc W}_\beta(\fg)$. As shown in
\cite{FR}, the latter arises in the limit of ${\cal W}_{q,t}({\fg})$
as $q \to 1$ and we impose the condition $t=q^\beta$.


In \cite{AFO}, we introduced the spaces of {\em deformed conformal
  blocks} of these two algebras for a curve ${\mc C}$ with the insertion
of vertex operators at finitely many marked points. We believe that
these spaces can only be defined in the case that the curve ${\mc C}$
is an infinite cylinder, or a plane, or a torus. But a plane can be
obtained from an infinite cylinder by taking the radius of the
cylinder to infinity, while the case of the torus should follow from
the case of the infinite cylinder by imposing periodic identifications
of marked points.

Suppose now that $^Lk \neq -^Lh^\vee$ (we will discuss the case
$^Lk=-^Lh^\vee$ in Section \ref{nonsimplylaced}). Let $q$ and $t$ in
${\cal W}_{q,t}({\fg})$ be given by the formulas \beq\label{zero}
q=\hbar^{-^L(k+h^{\vee})}, \qquad t=q^{m+\frac{1}{^L(k+h^\vee)}}. \eeq
Note that if we define $\beta$ by our relation \eqref{zero1}, then the
formulas in \eqref{zero} give \beq\label{first} \hbar = tq^{-m}, \qquad
t=q^\beta.\eeq

In \cite{AFO}, we proved, {\em if $\fg$ is simply-laced}, a
canonical isomorphism between the spaces of deformed conformal blocks
of $U_{\hbar}(\wh{^L\fg})_{^Lk}$ and ${\cal W}_{q,t}({\fg})$,
\beq\label{qLanglands} U_{\hbar}(\wh{^L\fg})_{^Lk}\on{-blocks} \qquad
\longleftrightarrow \qquad {\cal W}_{q,t}({\fg})\on{-blocks}, \eeq
where $q$ and $t$ are expressed in terms of $\hbar$ and $^Lk$ by
formula \eqref{zero}, and $\hbar$ and $^Lk$ are generic.

%
%

The case of {\em non-simply laced} $\fg$ is more subtle. We will
discuss it in Section \ref{nonsimplylaced} below.

The above identification of the deformed conformal blocks is what we
refer to as a {\em two-parameter Langlands correspondence} in the
title of the present paper (equivalently, the quantum $q$-Langlands
correspondence in the title of \cite{AFO}).

\begin{remark}
As we explained in \cite{AFO}, the physical setting for this
correspondence is a six-dimensional string theory, called the $(2,0)$
little string theory. This theory was introduced in \cite{SeibergN,
  Mm} and is a one-parameter deformation of the famous 6d $(2,0)$
superconformal theory (see e.g. \cite{Wittenl1}) -- the parameter is the characteristic size of the string. This deforms the relevant chiral algebras, associated to
$\ghat$ and ${\mc W}_\beta(\fg)$, into the quantum affine algebra and
the deformed ${\mc W}$-algebra, respectively. The little string interpretation of the deformed  ${\mc W}$-algebra is due to \cite{AH}.\qed
\end{remark}

In the case of an infinite cylinder ${\mc C}$, the space of
conformal blocks of $U(\Lfgh)_{^Lk}$ can be realized as the space of
solutions of the Knizhnik-Zamolodchikov (KZ) equations. We define the
space of deformed conformal blocks of the quantum affine algebra
$U_{\hbar}(\wh{^L\fg})_{^Lk}$ similarly, following \cite{FIR}, as the
space of solutions
of the quantum Knizhnik-Zamolodchikov (qKZ) equations. In both cases,
there is a particular fundamental solution of the equations which
comes from sewing chiral vertex operators. This solution is given by
formula \eqref{electric} in the case of deformed conformal blocks of
$U_{\hbar}(\Lfgh)_{^Lk}$.

On the other hand, the appropriate notion of the space of deformed
conformal blocks for the deformed ${\mc W}$-algebra ${\cal
  W}_{q,t}({\fg})$ was introduced by us in \cite{AFO}. These blocks
are defined using the free field realization of the vertex operators
of ${\mc W}_{q,t}({\fg})$ algebra in \eqref{magnetic} as well as the
screening operators, both defined in \cite{FR}. In addition, one has
to specify the contours of integration for the screening charges.

Thus, in \cite{AFO} we defined the space of deformed conformal blocks
of ${\mc W}_{q,t}({\fg})$ and we also introduced analogues of the qKZ
equations for these blocks. Our key insight was the {\em geometric
  interpretation} of these objects in terms of the (quantum)
$K$-theory of a Nakajima quiver variety $X$ \cite{OK}, whose quiver
diagram is based on the Dynkin diagram of ${\fg}$.

This paper is organized as follows. In Section \ref{s_statement}, we
define the deformed conformal blocks on both sides of the
two-parameter Langlands correspondence and state the main result,
Theorem \ref{t:one}. In Section \ref{s_Nakajima}, we explain the
geometry behind this correspondence; namely, that the deformed
conformal blocks can be realized as vertex functions in the
equivariant quantum K-theory of a Nakajima quiver variety. In Section
\ref{outline}, we give an outline of the proof of Theorem \ref{t:one}
following \cite{AFO}. Finally, in Section \ref{nonsimplylaced}, we
discuss a generalization to the non-simply laced case.

\section{Statement of the correspondence} 
\label{s_statement} 

Let $x$ be a coordinate on ${\cal C}\cong \Ct$. Fix a finite
collection of distinct points on ${\cal C}$, with coordinates $a_i$. 
We propose, and prove in the simply-laced case, 
a correspondence between the following two types of deformed conformal
blocks on ${\cal C}$.

\subsection{Electric side} 

\begin{figure}[!hbtp]
  \centering
   \includegraphics[scale=0.5]{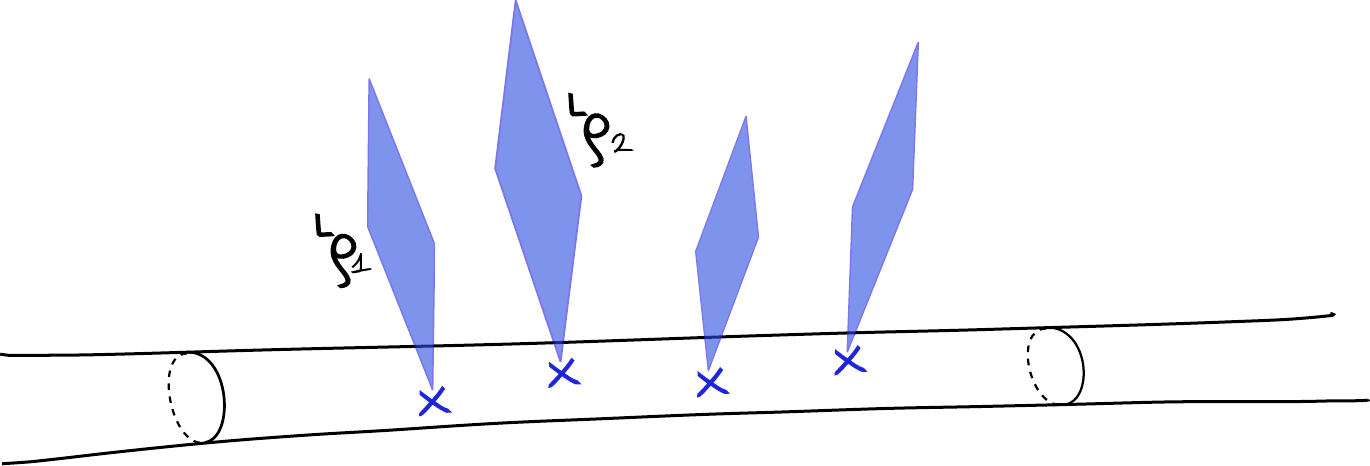}
 \caption{The cylinder ${\cal C}$ with the insertions of 
vertex operators corresponding to 
finite-dimensional $U_{\hbar}({{\Lfgh}})_{^Lk}$-modules $^L\!\rho_i$
at the points $a_i\in {\cal C}$. 
Boundary conditions at infinity are 
the highest weight vectors 
 $\langle \lambda'|$ and $|\lambda \rangle$.}
  \label{f_cylinder}
\end{figure}

On the \emph{electric} side, we consider the quantum affine algebra
$U_{\hbar}(\Lfgh)_{^Lk}$ blocks \cite{FIR}
\beq\label{electric}
\langle \lambda'| \;\prod_{i} { \Phi}_{{}^L\!\rho_i} (a_{i})  \;|\lambda \rangle
\eeq
where ${\Phi}_{{}^L\!\rho}(x)$ is a chiral vertex operator corresponding to a
finite-dimensional $U_{\hbar}({{\Lfgh}})_{^Lk}$-module $^L\!\rho$. 
The state $|\lambda\rangle
$ is the highest weight vector in a level $^Lk$ Verma module. 
Its weight $\lambda \in {^L\fh}^*$ is an element of the dual of the
Cartan subalgebra for $^L{\fg}$. This is illustrated in 
Figure \ref{f_cylinder}. 

In \cite{AFO}, we focused on the vertex operators corresponding to the
fundamental representations, with highest weights being the
fundamental weights $^Lw_a$. All others may be generated from them (by
fusion). The block \eqref{electric} takes values in
a weight subspace of
$$ 
\otimes_{ i} \, ({^L\!\rho}_i)= \otimes_a \,( ^L\!\rho_a)^{\otimes
  m_a} \,,
$$
namely, it has 
\begin{align}
  \textup{weight} & = \lambda' - \lambda \notag \\
 &= \label{weight}
\sum_a m_a\, ^Lw_a - \sum_a d_a \,^Le_a\,, 
\quad d_a \ge 0 \,. 
\end{align}
Here, the $^Le_a, a = 1,\ldots,{\rm rk}({\fg})$, are the simple
positive roots of $^L{\fg}$.

\subsection{Magnetic side}

On the \emph{magnetic} side, we consider matrix elements of
vertex and screening operators of ${\cal
  W}_{q,t}({\fg})$ algebra. They have the form
\beq\label{magnetic}
\langle \mu'| \;\prod_{i} {V}^{\vee}_ i(a_{i}) \prod_a \Bigl({Q}_a^{\vee}\Bigr)^{d_a} \;|\mu \rangle.
\eeq
where $V^{\vee}_a(x)$ and $Q^{\vee}_a$ are the vertex and the
screening charge operators, respectively, defined in \cite{FR}. They
are labeled by the coroots and
coweights of ${\fg}$, respectively. The screening charge operators are
defined as integrals of screening current vertex operators $Q^{\vee}_a
= \int dx \; S_a^{\vee}(x)$, so \eqref{magnetic} is in fact an
integral formula, in which contours of integration must be specified
(and this is a non-trivial problem).

The coweights of ${\fg}$ labeling the $V^{\vee}_a(x)$ are the highest
weights of the fundamental representations of $^L{\fg}$.  The operator
${V}^{\vee}_ i(a_{i})$, inserted at a point on ${\cal C}$ with the
coordinate $a_i$, is associated to the same representation of
$^L{\fg}$ as the corresponding vertex operator in
\eqref{electric}. The state $|\mu\rangle$, labeled by an element $\mu
\in {\fh}$ of the Cartan subalgebra of ${\fg}$, generates an
irreducible Fock representation of the ${\cal W}_{q,t}({\fg})$ algebra
\cite{FR}. The (co)weights $\mu$ and $\mu'$ are
  determined by $\lambda$ and $\lambda'$ (the exact formula depends on
  the chosen normalization).


\subsection{Main result}

The main result of \cite{AFO} is the following theorem:

\begin{Theorem}\label{t:one} Let ${\fg}$ be a simply-laced simple Lie algebra.
The deformed conformal blocks of $U_{\hbar}({{\Lfgh}})_{^Lk}$ in
\eqref{electric} and the deformed conformal blocks of ${\mc
  W}_{q,t}({\fg})$ in \eqref{magnetic}, whose parameters are generic
and are related by equation \eqref{zero}, are canonically identified by
the formula \beq\begin{aligned}\label{linear} \textup{specific
  covector} \times \textup{ $U_{\hbar}(\Lfgh)_{^Lk}$-block} \quad \quad
\leftrightarrow \quad \quad \textup{${\cal W}_{q,t}(\fg)$-block}.
\end{aligned}
\eeq
  \end{Theorem}

The covector in \eqref{linear}, as well other ingredients of Theorem
\ref{t:one}, are explained in \cite{AFO} in geometric terms; namely, in
terms of the (quantum) K-theory of a Nakajima quiver variety $X$ (see
the next section). Specifically, the covector in question corresponds
to the insertion of the identity $\cO_X\in K_\bT(X)$ in a certain
enumerative problem.

\section{{Geometry behind the correspondence}}
\label{s_Nakajima}

The central ingredient of our proof is that for Lie algebras of
simply-laced type, i.e.
$$
^L{\fg} ={\fg},
$$
we can realize the deformed conformal blocks \eqref{electric} and
\eqref{magnetic} as {\it vertex functions} in the equivariant 
quantum K-theory of a
certain holomorphic symplectic variety $X$. The variety $X$ is the
{\em Nakajima quiver variety} with 
$$
\textup{quiver ${\scQ}$} = 
\textup{Dynkin diagram of ${\fg}$} \,.
$$

\subsection{Nakajima quiver variety}

The Nakajima quiver variety $X$ is a hyper-K\"ahler  quotient
(or a holomorphic symplectic reduction) 
\beq\label{NV}
X=T^*{{\rm Rep }\,\scQ} \rd G_{\scQ}\,, 
\eeq
where 
\beq\label{repQ}
{{\rm Rep }\,\scQ}=\oplus_{a \rightarrow b} \,{\rm Hom}(V_a, V_b) \oplus_{a} {\rm Hom}(V_a, W_a)
\eeq
and 
\beq\label{gaugeC}
G_{\scQ} =\prod_{a} GL(V_a), \qquad G_{W} = \prod_{a} GL(W_a) \,. 
\eeq
The arrows in \eqref{repQ} are the arrows of the quiver. The
dimensions of the vector spaces $V_a$ and $W_a$ correspond 
as follows 
$$
\dim V_a = d_a\,, \quad \dim W_a = m_a
$$
to the weight space data in \eqref{weight}. 

The quotient in \eqref{NV} involves a geometric invariant theory (GIT)
quotient, which depends on a choice of stability conditions.
As a result, vertex functions also depend on a stability
condition. This stability condition makes them analytic in a certain 
region of the K\"ahler moduli space of $X$. The transition matrix between
the vertex functions and the deformed conformal blocks also 
depend on the stability condition.

\subsection{Equivariant variables}

Most variables in \eqref{electric} and \eqref{magnetic} 
become \emph{equivariant} in their
geometric interpretation. We have 
$$
G_W\times \Ct_\hbar \subset \Aut(X)
$$
where $\Ct_\hbar$ rescales the cotangent directions in \eqref{NV} with
the weight $\hbar^{-1}$. This gives the symplectic form on $X$ weight
$\hbar$ under $\Ct_\hbar$. We fix a maximal torus 
$A\subset G_W$ and denote 
$$
T = A \times \Ct_\hbar \,. 
$$
The coordinates $a_i$ of $A$ are the positions at which the vertex
operators are inserted in \eqref{electric} and \eqref{magnetic}, 
while $\hbar$ is the quantum group deformation parameter in 
$U_{\hbar}(\Lfgh)_{^Lk}$. 

A multiplicative group $\Ct_q$ acts on quasimaps
$\bP^1\dasharrow X$ by automorphisms
of the domain $\bP^1$. The coordinate $q\in \Ct_q$ is the
$q$-difference parameter from the title of the paper.

In \cite{Nakajima}, Nakajima identified $K_T(X)$ with a space of 
weight \eqref{weight} 
in a $U_{\hbar}(\Lfgh)_{^Lk}$-module.

\subsection{Vertex functions}    \label{qKZ}

The basic object of the theory of \cite{OK} is the {\em vertex
  function} ${\bf V}$. The vertex
  function is an equivariant $K$-theoretic count of quasimaps from ${\mathbb C}$ to $X$
 of all possible degrees.   It is an analog of Givental's ${\bf
  I}$-function. 
The variables $z$ in this generating function
are called \emph{K\"ahler parameters}. They are related to the choice 
of the Fock vacuum  $|\lambda\rangle$ in \eqref{electric} and $|\mu\rangle$ in
\eqref{magnetic}.  
For its definition and basic properties, see \cite{AFO}, Section 3.2.

A key geometric property of the vertex functions are the 
$q$-difference
equations that they satisfy, as functions of both equivariant and 
K\"ahler variables (see \cite{OK}, Section 8, for an introduction).
In particular, the $q$-difference equations in the variables $a_i$ were 
identified in Section 10 of \cite{OK} with the qKZ equations
\cite{FIR}. In \cite{FIR}, these were introduced as the $q$-difference
equations that determine the deformations of conformal blocks
corresponding to ${\Lfgh}$ in \eqref{electric}. 

More precisely, the fundamental solutions of the qKZ equations are
vertex functions counting maps from
$
{\mathbb C}^{\times}$ to  $X$
together with relative insertions at $0 \in \C$
\cite{OK}.

\section{Outline of the proof of the main theorem}    \label{outline}

The proof of Theorem \ref{t:one} follows by establishing connections
between \eqref{electric}, \eqref{magnetic}, and the vertex
functions. Here's an outline.

\subsection{Electric side}

On the electric side, i.e. the
$U_{\hbar}(\Lfgh)_{^Lk}$-algebra side, we have a
  characterization of deformed conformal blocks in \eqref{electric} as
solutions of the qKZ equations that they satisfy. 
Vector vertex functions provide a different basis of solutions of 
the same qKZ equations. The difference manifests itself through
difference analytic dependence on the equivariant variables $a_i$ and 
the K\"ahler variables $z$.
Like any two bases of meromorphic solutions to the same difference 
equations, the vector vertex functions and the $U_{\hbar}(\Lfgh)_{^Lk}$-blocks 
are connected by a $q$-periodic transition matrix. This $q$-periodic transition matrix
may be called the \emph{pole
  subtraction matrix} because it cancels unwanted poles 
in one set of variables at the expense of introducing poles in another 
set of variables, see \cite{ese} for a detailed discussion.

This pole subtraction matrix was identified geometrically in \cite{ese} 
as the elliptic cohomology version of the stable envelopes of the Nakajima
variety $X$. Stable envelopes in equivariant 
cohomology were introduced in \cite{MO}. They
are the main geometric input in the construction of quantum group
actions suggested there, see Section 9 of \cite{OK} for an overview.
This notion has a natural lift to equivariant 
K-theory, derived categories of coherent sheaves, and, as shown in 
\cite{ese}, also to the equivariant elliptic 
cohomology. 

In parallel to cohomology and K-theory, elliptic stable envelopes
produce an action of a quantum group, namely an elliptic quantum
group. The analysis of \cite{ese} relates the monodromy of the qKZ
equations to the braiding for this elliptic quantum group. First
steps towards such identification were taken already in \cite{FIR},
with many subsequent developments, as discussed in \cite{ese}.

In the enumerative problem, elliptic stable envelopes are inserted via
the the evaluation map at infinity of $\C^\times$, away from the point
$0$ where the relative conditions have been inserted. They appear as
elliptic functions multiplying the measure of integration (see
\cite{AFO}, Section 2.2). In either interpretation, they map vector
vertex functions to $U_{\hbar}(\Lfgh)_{^Lk}$-blocks.

\subsection{Magnetic side}

On the magnetic, i.e. $\mathcal{W}$-algebra side, we prove in
\cite{AFO}, Theorem 3.1 that the vertex functions ${\bf V}$ of $X$,
counting quasimaps \beq\label{scalar}\C \dasharrow X, \eeq equal the
integrals \eqref{magnetic} for a specific choices of contours of
integration. The integral formulas for vertex
functions of $X$ arise as follows.

K-theoretic computations on a GIT-quotient by a
reductive group $G$ may be expressed as $G$-invariants in a 
$G$-equivariant computation on the prequotient. The projection onto 
$G$-invariants may be recast, by the Weyl integration formula, as an
integral over a suitable cycle in a maximal torus in $G$.
In \cite{AFO}, we show that
for K-theoretic computations on the moduli spaces of quasimaps to a 
GIT-quotient, there are similar integral formulas.

\subsection{The match of deformed conformal blocks} 

To establish the match, it suffices
to recognize in these formulas the integral formulas of \cite{FR} for
the free field correlators of ${\cal W}_{q,t}({\fg})$.

The same dichotomy arises in the discussion of the magnetic deformed
conformal blocks. Vertex functions are analytic as $z\to 0$, while the
natural requirement for the ${\cal W}_{q,t}({\fg})$-blocks is to be
analytic in regions of the form
\begin{equation}
|a_5| \gg |a_1| \gg |a_3| \gg \dots,\label{a_region}
\end{equation}
Importantly, the very
same elliptic stable envelopes transform the $z$-series into functions
with the right analyticity in the $a$-variables. The geometry of the
correspondence is tautologically the same, as the insertion of the
elliptic stable envelope happens at infinity, away from the point $0$,
which distinguishes vertex functions from their vector analogs.  In
integral formulas, stable envelopes appear as elliptic functions
multiplying the measure of integration.

Informally, vertex functions are a special case of the vector vertex
functions, namely the one corresponding to no insertion at $0$. Since
the moduli spaces in questions are not really identical, the correct
technical way to see this it is via the degeneration formula as in
\cite{AFO}, Section 4.1. Then applying elliptic stable envelopes to
both sides, we obtain the statement of Theorem \ref{t:one}. See
\cite{AFO} for more details.

\section{Non-simply laced case}    \label{nonsimplylaced}

If $\fg$ is non-simply laced, the statement of Theorem \ref{t:one}
should be modified in a non-trivial way. On the magnetic side of the
isomorphism \eqref{linear} we still have the deformed conformal blocks
of ${\mc W}_{q,t}({\fg})$ (given by formula \eqref{magnetic}) but what
should appear on the electric side is {\em not} the deformed conformal
blocks of $U_{\hbar}({{\Lfgh}})_{^Lk}$ (given by formula
\eqref{electric}).

To see that, consider the critical level limit of
$U_{\hbar}({{\Lfgh}})_{^Lk}$, i.e. the limit in which its level $^Lk$
tends to $-^Lh^\vee$. In this limit, the first formula in \eqref{zero}
yields $q=1$, and the first formula in \eqref{first} yields $\hbar=t$.

It is known that in this limit, the solutions of the qKZ equations of
$U_{\hbar}({{\Lfgh}})_{^Lk}$ give rise to the Bethe eigenvectors of
the XXZ-type model associated to the quantum affine algebra
$U_{\hbar}({{\Lfgh}})$. Moreover, the difference operators of the qKZ
equations become in this limit the transfer-matrices, i.e. the
Hamiltonians, of this model (see the Appendix of \cite{FHR} for
details). Thus, the system of qKZ equations can be viewed as a
deformation of this model.

If $\fg$ is simply-laced, then ${\mc W}_{1,t}(\fg)$ is isomorphic to
the algebra of transfer-matrices of the XXZ-type model associated to
$U_t(\ghat) = U_\hbar({{\Lfgh}})$, so everything is consistent.
However, if $\fg$ is not simply-laced, the algebra ${\mc
  W}_{1,t}(\ghat)$ is not isomorphic to the algebra of
transfer-matrices of the XXZ-type model associated to
$U_t(\Lfgh)$. Rather, ${\mc W}_{1,t}(\ghat)$ is isomorphic
to the algebra of transfer-matrices of a different model, the {\em
  folded quantum integrable model} associated to $\fg$, which was
defined in \cite{FHR} by ``folding'' the XXZ-type model associated to
$U_t(\wh{\fg'})$, where $\fg'$ is the simply-laced simple Lie
algebra equipped with an automorphism of order 2 or 3 whose invariant
Lie subalgebra is $\fg$.

This suggests that for non-simply laced $\fg$ the two-parameter
Langlands correspondence is more subtle. Namely, it follows from the
preceding paragraph that the $q$-difference equations appearing in
this correspondence are {\em not} the usual qKZ equations associated
to $U_\hbar(\Lfgh)_{^Lk}$ (if they were, then in the critical level
limit we would recover the eigenvectors of the XXZ-type model
associated to $U_\hbar(\Lfgh) = U_t(\Lfgh)$, but this would be
inconsistent with the limit on the other side of the correspondence,
which yields ${\mc W}_{1,t}(\fg)$). Rather, it must be a {\em new
  $q$-deformation of the KZ equation of $ \Lfgh_{^Lk}$}, yet to be
defined.

We are going to address this issue in our forthcoming paper \cite{AFO1}.


\begin{thebibliography}
\raggedright

\bibitem{AFO} M.~Aganagic, E.~Frenkel, and A.~Okounkov,
  ``Quantum $q$-Langlands correspondence,''
  Trans. Moscow Math. Soc. {\bf 79} (2018) 1-83
  [arXiv:1701.03146].

\bibitem{AFO1} M.~Aganagic, E.~Frenkel, and A.~Okounkov,
  ``Quantum $q$-Langlands correspondence in the non-simply laced case,''
to appear.

\bibitem{AH}
M.~Aganagic and N.~Haouzi,
``ADE Little String Theory on a Riemann Surface (and Triality),''
[arXiv:1506.04183].
  
\bibitem{ese} 
  M.~Aganagic and A.~Okounkov,
  ``Elliptic stable envelope,''
  [arXiv:1604.00423].

\bibitem{Bethe} M.~Aganagic and A.~ Okounkov, ``Quasimap counts and
     Bethe eigenfunctions,''  
Mosc.\ Math. J.\ \textbf{17} (2017) 565-600. 

\bibitem{BD}
A.~Beilinson and V.~Drinfeld, ``Quantization of Hitchin's Integrable
System and Hecke Eigensheaves,'' Preprint at
\url{https://www.math.uchicago.edu/~arinkin/langlands}
 
\bibitem{BO} M. Bershadsky and H. Ooguri, ``Hidden $SL(n)$ symmetry in
  conformal field theory,'' Commun. Math. Phys. {\bf 126} (1989) 49-83.
  
\bibitem{Drinfeld} V. Drinfeld, ``Hopf algebras and the quantum
  Yang-Baxter equation,'' Dokl. Soviet Math. {\bf 32} (1985) 254-258.

\bibitem{EFK} P. Etingof, E. Frenkel, and D. Kazhdan, ``A general
  framework for the analytic Langlands
  correspondence,'' Pure and Applied
Mathematics Quarterly {\bf 20} (2024) 307-426
[arXiv:2311.03743].
  
   \bibitem{FF:ds} B. Feigin and E. Frenkel, ``Quantization of the
  Drinfeld-Sokolov reduction,'' Phys. Lett. {\bf B246} (1990) 75-81.

\bibitem{FF} B. Feigin and E. Frenkel, ``Affine Kac-Moody algebras at
  the critical level and Gelfand-Dikii algebras,''
  Int. J. Mod. Phys. {\bf A7}, Suppl. 1A (1992) 197-215 (reprinted in
  ``Infinite analysis'' (Kyoto, 1991), Adv. Ser. Math. Phys. 16, World
  Scientific, 1992).

\bibitem{F:book} E. Frenkel, ``Langlands Correspondence for Loop
  Groups,''
  Cambridge Studies in Advanced Mathematics {\bf 103}, Cambridge
  University Press, 2007.

\bibitem{Frenkel}
{E.~{Frenkel}}, ``Lectures on the Langlands Program and Conformal
Field Theory,'' in Frontiers in Number Theory, Physics and Geometry
    II, eds. P. Cartier, e.a., pp. 387-536, Springer, 2007
[arXiv:hep-th/0512172].

\bibitem{FB} E.~Frenkel and D.~Ben-Zvi, ``Vertex Algebras and
Algebraic Curves,'' Mathematical Surveys and
Monographs {\bf 88}, Second Edition, AMS, 2004.

  \bibitem{FHR} E.~Frenkel, D.~Hernandez, and N.~{Reshetikhin},
    ``Folded quantum integrable models and deformed ${\mathcal
    W}$-algebras,''
    Letters in Math. Physics {\bf 112} (2022) 80 [arXiv:2110.14600].
    
\bibitem{FR} E.~Frenkel and N.~{Reshetikhin}, ``Deformations of
    W-algebras associated to simple Lie algebras,''
    Comm. Math. Phys. {\bf 197}
(1998) 1-32 [arXiv:q-alg/9708006].

  \bibitem{FIR}
I.B.~Frenkel and N.Yu.~Reshetikhin, ``Quantum affine algebras and
holonomic difference equations,'' Comm. Math. Phys. {\bf 146} (1992) 1-60.

\bibitem{Gaitsgory} D. Gaitsgory, S. Raskin, e.a., ``Proof of the geometric
  Langlands conjecture,'' \newline Preprints at
  \url{https://people.mpim-bonn.mpg.de/gaitsgde/GLC/}

\bibitem{Jimbo} M. Jimbo, ``A $q$-difference analogue of
  $\mathcal{U}({\mathfrak{g}})$ and Yang-Baxter equation,''
  Lett. Math. Phys. {\bf 10} (1985) 63-69.

\bibitem{KW}
  A.~Kapustin and E.~Witten,
  ``Electric-Magnetic Duality And The Geometric Langlands Program,''
  Commun.\ Num.\ Theor.\ Phys.\  {\bf 1} (2007) 1-236
  [arXiv:hep-th/0604151].

\bibitem{Mm} 
  A.~Losev, G.~W.~Moore, and S.~L.~Shatashvili,
  ``M \& m's,''
  Nucl.\ Phys.\ {\bf B522} (1998) 105-124
  [arXiv:hep-th/9707250].

\bibitem{MO} 
  D.~Maulik and A.~Okounkov,
  ``Quantum Groups and Quantum Cohomology,''
  Ast\'erisque {\bf 408} (2019) [arXiv:1211.1287].

\bibitem{Nakajima}
H.~Nakajima, ``Quiver varieties and finite-dimensional representations
of quantum affine algebras,''
{J. Amer. Math. Soc.}, {\bf 14}, (2001) 145-238.

 \bibitem{OK}
A.~Okounkov, ``Lectures on K-theoretic computations in enumerative
geometry,'' Geometry of moduli spaces and representation theory, 251-380, 
IAS/Park City Math.\ Ser., 24, AMS 2017. 

\bibitem{SeibergN} 
  N.~Seiberg,
  ``New theories in six-dimensions and matrix description of M theory on $T^5$ and $T^5/Z_2$,''
  Phys.\ Lett.\ {\bf B408} (1997) 98-104.
  [arXiv:hep-th/9705221].

\bibitem{Wittenl1} 
  E.~Witten,
  ``Geometric Langlands From Six Dimensions,'' [arXiv:0905.2720].

\end{thebibliography}
\end{document}